\documentclass{amsart}
\usepackage{amsfonts}
\usepackage{amsmath,amscd}
\usepackage{amssymb}
\usepackage{amsthm}
\usepackage{newlfont}

 \newtheorem{thm}{Theorem}[section]
 \newtheorem{cor}[thm]{Corollary}
 \newtheorem{lem}[thm]{Lemma}
 \newtheorem{prop}[thm]{Proposition}
 \theoremstyle{definition}
 
 \theoremstyle{remark}
 \newtheorem{rem}[thm]{Remark}

 \numberwithin{equation}{section}

\begin{document}

\title[On an argument of J.--F. Cardoso dealing with ...]
 {On an argument of J.--F. Cardoso dealing with perturbations of joint diagonalizers}

\author[F.G.  Russo]{Francesco G. Russo}
\address{\newline
Department of Mathematics\newline
University of Palermo\newline
via Archirafi 34, 90123, Palermo, Italy \newline 
URL: www.fgrusso.altervista.org \newline
E-mail: francescog.russo@yahoo.com}



\date{\today}

\begin{abstract}
B. Afsari has recently proposed a new approach to the matrix joint diagonalization, introduced by J.--F. Cardoso in 1994,  in order to investigate the independent component analysis and the blind signal processing in a wider prospective.  Delicate notions of linear algebra and differential geometry are involved in the works of B. Afsari and the present paper continues such a line of research, focusing on a theoretical condition which has significant consequences in the numerical applications.
\end{abstract}

\maketitle

\section{EJD problem and relative generalizations}
Several branches of the  digital signal processing involve linear algebra, since the empirical data are often organized in suitable  matrices,  which describe the properties of  some signals with respect to others, rich of noise and then less significative. Roughly speaking, this  idea is behind the construction of the filters, which play a fundamental role in most of the mathematical models and of the simulations, currently used in the applied sciences. Therefore it is not suprising whether some aspects of  music and of acoustics are treated in \cite{r2,r3,r4} in terms of  linear systems. A classical reference remains \cite{pm} for the description of the common techniques of the digital signal processing. Among these, the blind signal processing and the indipendent component analysis, more properly treated in \cite{a1,a2,a3,a4,bbm,c1,c2,hls}, have been largely investigated in the last years for various motivations and many algorithms have been developped.

A significant contribution on the blind signal processing and on the indipendent component analysis is \cite{a4}, where the original ideas of J.--F. Cardoso \cite{c1,c2} have been generalized, introducing the problems of Exact Joint Diagonalization (EJD), Orthogonal Joint Diagonalization (OJD) and Non--Orthogonal Joint Diagonalization (NOJD). 
 \cite[Equation 3.1]{a4} shows that the OJD poblem, originally introduced in \cite{c1}, may be reduced to a cost function for which conditions of existence and unicity (in the sense of \cite[Theorem 2.4]{a4}) should be examined. A similar reduction may be done for the NOJD problem, as described in \cite[\S 3.2.2 and \S 3.2.3]{a4}, looking for conditions as in \cite[Theorem 4.1]{a4}.

Therefore a problem of applied science becomes a problem of optimization, which can be approached with the methods of the linear algebra,  looking for suitable matrices which satisfy  a prescribed equation. Analogies are very common in literature and \cite{hm,sj}   illustrate the relations among riemannian geometry,   topology,  theory of compact Lie groups and numerical analysis, when we want to solve some problems of optimization as those which appears in the blind signal processing and in the indipendent component analysis.  In spite of the  algorithms which can be written with sophisticated computer programs, we will concentrate only on a theoretical aspect of an equation of matrices.

Section 2 recalls some feedback and notations from \cite{hm,hls,r,sj}, which will allow us to give short proofs of \cite[Lemmas 1,2, 3] {c1},  fundamental in  \cite{c1} for the EJD problem. Section 3 deals with the main result, which is a generalization of \cite[Proposition 1]{c1}. Final considerations and open problems are placed at the end.

\section{Joint diagonalizers revisited}
We begin with some well--known notions on groups of matrices, referring  to \cite{hm,r,sj} for the terminology.

In the Hilbert space $\mathbb{C}^n$ of  dimension $n\ge 1$ over the field $\mathbb{C}$ of the complex numbers we have the standard
scalar product $(\cdot \ | \ \cdot):  (x,y)=((x_1,\ldots x_n),(y_1,\ldots y_n))\in \mathbb{C}^{n\times n}
\mapsto (x,y)=x_1\overline{y}_1+\ldots +x_n\overline{y}_n\in [0,+\infty[$. The set $\mathrm{GL}(n,\mathbb{C})$
of all non--singular $n\times n$ matrices with coefficients in $\mathbb{C}$ is a non--abelian group, called \textit{general linear group  of dimension n over $\mathbb{C}$}.
The same is true, when we replace $\mathbb{C}$ with the field $\mathbb{R}$ of the real numbers and consider the
corresponding standard scalar product. Consequently, the following groups are defined
$\mathrm{GL}(n,\mathbb{C})=\{B\in \mathbb{C}^{n\times n} \ | \ \mathrm{det}(B)\not=0\}$  and
$\mathrm{GL}(n,\mathbb{R})=\{B\in \mathbb{R}^{n\times n} \ | \ \mathrm{det}(B)\not=0\}$, where $\mathrm{det}(B)$ denotes the  determinant of $B$ in the usual sense. 

The set $Z(\mathrm{GL}(n,\mathbb{C}))=\{A,B \in
\mathrm{GL}(n,\mathbb{C}) \ | \ AB=BA\}$ is the \textit{center of} $\mathrm{GL}(n,\mathbb{C})$ and it is straightforward to see that 
$Z(\mathrm{GL}(n,\mathbb{C}))=\{aI\in \mathrm{GL}(n,\mathbb{C}) \ | \ a\in
\mathbb{C}-\{0\} \}.$
The quotient group
$\mathrm{PGL}(n,\mathbb{C})=\mathrm{GL}(n,\mathbb{C})/Z(\mathrm{GL}(n,\mathbb{C}))$
is the \textit{projective linear group of dimension n over $\mathbb{C}$}. By analogy,
$\mathrm{PGL}(n,\mathbb{R})=\mathrm{GL}(n,\mathbb{R})/Z(\mathrm{GL}(n,\mathbb{R})).$

On another hand,  
$\mathrm{SL}(n,\mathbb{C})=\{B\in \mathbb{C}^{n\times n} \ | \ \mathrm{det}(B)=1\}$ (resp.    
$\mathrm{SL}(n,\mathbb{R})=\{B\in \mathbb{R}^{n\times n} \ | \ \mathrm{det}(B)=1\}$)
is a subgroup of $\mathrm{GL}(n,\mathbb{C})$ (resp.  $\mathrm{GL}(n,\mathbb{R})$), called \textit{special linear group of dimension n over $\mathbb{C}$} (resp. $\mathbb{R}$). This time 
$Z(\mathrm{SL}(n,\mathbb{C}))=\{aI\in
\mathrm{GL}(n,\mathbb{C}) \ | \ a^n=1,  \ a\in \mathbb{C}-\{0\} \}$ and
$\mathrm{PSL}(n,\mathbb{R})=\mathrm{SL}(n,\mathbb{C})/Z(\mathrm{SL}(n,\mathbb{C}))$ is called \textit{projective special linear group of dimension n over $\mathbb{C}$}; similarly for $\mathrm{PSL}(n,\mathbb{R})=\mathrm{SL}(n,\mathbb{R})/Z(\mathrm{SL}(n,\mathbb{R}))$.

$\mathrm{U}(n)=\{B\in \mathrm{GL}(n,\mathbb{C}) \ | \ (Bx,Bx)=(x,x), x\in \mathbb{C}^n\}=\{B\in
\mathrm{GL}(n,\mathbb{C}) \ | \ B^{-1}=B^H\},$ 
where $B^H=\overline{B}^*$ denotes the conjugate
transposed of $B$, is  a subgroup of $\mathrm{GL}(n,\mathbb{C})$, called  \textit{unitary group of dimension
n over $\mathbb{C}$}, and the intersection $\mathrm{U}(n)\cap \mathrm{SL}(n,\mathbb{C})=\mathrm{SU}(n)$ is still
a subgroup of $\mathrm{GL}(n,\mathbb{C})$, called \textit{special unitary group of dimension n over
$\mathbb{C}$}. By analogy, 
$\mathrm{O}(n)=\{B\in \mathrm{GL}(n,\mathbb{R}) \ | \ (Bx,Bx)=(x,x),
x\in \mathbb{R}^n\}=\{B\in \mathrm{GL}(n,\mathbb{R}) \ | \ B^{-1}=B^*\},$
is the \textit{orthogonal group of dimension n over $\mathbb{R}$} and $\mathrm{SO}(n)=\mathrm{O}(n)\cap
\mathrm{SL}(n,\mathbb{R})$ is the \textit{special orthogonal group of dimension n over $\mathbb{R}$}.
Furthermore, $B\in \mathrm{GL}(n,\mathbb{C})$ is called \textit{hermitian}, if $B^H=B$, and
\textit{anti--hermitian}, if $B^H=-B$.

$\mathrm{GL}(n,\mathbb{C})$, endowed with the topology induced by
the operator norm $\| \ \| : A \in \mathrm{GL}(n,\mathbb{C}) \mapsto \|A\|=\mathrm{sup}\{\|Ax\| \ | \ \|x\|\le1
\}\in [0,+\infty[$ for all $x \in \mathbb{C}^n$, is a topological group. The same is true for
$\mathrm{GL}(n,\mathbb{R})$, mutatis mutandis. In particular, all groups $\mathrm{O}(n)$, $\mathrm{SO}(n)$, $\mathrm{U}(n)$, $\mathrm{SU}(n)$ are compact groups for $n=1,2,\ldots$.

All these notions allow us to reformulate more properly the crucial points in \cite{a4,c1} and to simplify some technical arguments.

If $A\in \mathrm{GL}(n,\mathbb{C})$, then we follow \cite{c1}, defining \textit{the off of $A$} as the map
\begin{equation}\label{e:1}
\mathrm{off}: A \in  \mathrm{GL}(n,\mathbb{C}) \mapsto  \mathrm{off}(A)={\underset{i\not=j} {\underset{ i,j \in
\{1,2,\ldots,n\}} \sum} }|a_{ij}|^2 \in [0,+\infty[,
\end{equation}
where $a_{ij}$ is the $(i,j)$--th entry of $A$. More generally, if $V \in \mathrm{U}(n)$ and
\begin{equation}\mathcal{M}=(M_1,M_2,\ldots, M_k,\ldots, M_m)\end{equation} is an $m$--tuple of matrices in $\mathrm{GL}(n,\mathbb{C})^m$, then
we may consider the map
\begin{equation}\label{e:2}
\mathcal{Y}: (V,\mathcal{M}) \in \mathrm{U}(n) \times \mathrm{GL}(n,\mathbb{C})^m \mapsto
\mathcal{Y}(V,\mathcal{M})=\sum^m_{k=1} \mathrm{off}(V^HM_kV)\in [0,+\infty[
\end{equation}
and $V$ is called \textit{unitary minimizer} (or \textit{joint diagonalizer}) of $\mathcal{Y}$ with respect to $\mathcal{M}$, whenever $\mathcal{Y}(V,\mathcal{M})=0$. In case of connected compact groups, existence and uniqueness of non--trivial solutions of the equation $\mathcal{Y}(V,\mathcal{M})=0$ are ensured by classical theorems of calculus on smooth manifolds (for instance, Weiertrass theorems, see \cite{hm}). Unfortunately, in general it is very hard to decide for which $(V,\mathcal{M})$ the equation $\mathcal{Y}(V,\mathcal{M})=0$ is fulfilled (namely, this is the EJD problem in \cite{a3,a4}).

Now we consider, instead of $\mathcal{M}$, the $m$-tuple of  $\mathrm{GL}(n,\mathbb{C})^m$
\begin{equation}\label{e:2bis}\mathcal{M}_0=(UD_1U^H, \ldots, UD_kU^H,\ldots, UD_mU^H)
\end{equation}
where $U\in \mathrm{U}(n)$ and $D_k$ is the diagonal matrix of $\mathbb{C}^{n\times n}$ with diagonal entries
$d_1(k),\ldots,d_n(k)$, that is,
\begin{displaymath}\label{displ:1}
D_1=\left(
\begin{array}{cccccccc}
d_1(1) & 0 & \ldots & 0 \\
0 & d_2(1)  & \ldots & 0 \\
\ldots & \ldots \\
0 & 0 & \ldots &  d_n(1)
\end{array}\right),
D_2=\left(
\begin{array}{cccccccc}
d_1(2) & 0 & \ldots & 0 \\
0 & d_2(2)  & \ldots & 0 \\
\ldots & \ldots \\
0 & 0 & \ldots &  d_n(2)
\end{array}\right),
\ldots \end{displaymath}
\begin{displaymath}\label{displ:2}
\ldots \ldots \ldots \ldots \ldots \ldots,  \ \ D_k=\left(
\begin{array}{cccccccc}
d_1(k) & 0 & \ldots & 0 \\
0 & d_2(k)  & \ldots & 0 \\
\ldots & \ldots \\
0 & 0 & \ldots &  d_n(k)
\end{array}\right), \ \ldots \ldots \ldots \ldots \ldots \ldots
\end{displaymath}
\begin{displaymath}
\ldots \ldots \ldots \ldots \ldots \ldots \ldots \ldots, \ D_m=\left(
\begin{array}{cccccccc}
d_1(m) & 0 & \ldots & 0 \\
0 & d_2(m)  & \ldots & 0 \\
\ldots & \ldots \\
0 & 0 & \ldots &  d_n(m)
\end{array}\right).
\end{displaymath}
We have \begin{equation}\label{trivial}\mathcal{Y}(U,\mathcal{M}_0)=\sum^m_{k=1} \mathrm{off}(U^H(UD_kU^H)U)=\sum^m_{k=1}
\mathrm{off}(U^{-1}UD_kU^{-1}U)=\sum^m_{k=1} \mathrm{off}(D_k)=0,\end{equation} and, of course, the map $\mathcal{Y}$ is
minimized. In this case it is immediate to see that $U$ is a joint diagonalizer for $\mathcal{Y}$ with respect to $\mathcal{M}_0$.

A more general situation can be encountered when we consider
\begin{equation}\label{e:3}
\mathcal{M}_\lambda=(UD_1U^H + \lambda R_1, \ldots, UD_kU^H + \lambda R_k, \ldots, UD_mU^H + \lambda R_m),
\end{equation}
where $\lambda \in \mathbb{R}$ and $R_k\in \mathrm{GL}(n,\mathbb{C})$ (for $k=1,2,\ldots,m$).  Of course, for $\lambda=0$ (\ref{e:3})
becomes (\ref{e:2bis}). Therefore it is very interesting to see (under the point of view of the EJC problem) if there are $\lambda\not=0$ for which  $(U,\mathcal{M}_\lambda) \in \mathrm{U}(n) \times \mathrm{GL}(n,\mathbb{C})^m$  satisfies $\mathcal{Y}(U,\mathcal{M}_\lambda)=0$. The answer is positive in the sense of  \cite[Proposition 1]{c1}. On another hand, small values of $\lambda$, close to 0, allow us to 
study  the perturbations in a neighborhood of the minimum value reached by $\mathcal{Y}$ on $(U,\mathcal{M}_0)$. Some consequences of this condition of equilibrium 
can be found in \cite{a1,a2,a3,a4,bbm,c2}, where it is emphasized the application to the blind beamforming of non--gaussian signals.

In order to state more properly \cite[Proposition 1]{c1}, we recall that  $e_i$ denotes the $n\times 1$ row vector with 1 in the $i$--th position and 0 elsewhere, that is, $e_1=(1,0,0,\ldots,0)$, $e_2=(0,1,0,\ldots,0)$, $\ldots$, $e_i=(0,\ldots,0,1,0,\ldots,0)$, $\ldots$, $e_n=(0,0,\ldots,1)$. $e^*_i$ denotes the $1\times n$ column vector with 1 in the $i$--th position and 0 elsewhere. This leads to the matrix 
\begin{equation}
\left(
\begin{array}{cccccccc}
0 & 0 & 0 &\ldots & 0 \\
0 & 0  & 1 &\ldots & 0 \\
\ldots & \ldots & \ldots & \ldots\\
0 & 0 & \ldots & \ldots &  0
\end{array}\right) 
\end{equation}
whose entries $e_{ij}$ are 1 in the $(i,j)$--th position and 0 elsewhere.
More generally, $a_i$ is the $i$--th row and $a^*_i$ the $i$--th column of $A \in \mathrm{GL}(n,\mathbb{C})$, whose entries are  $a_{ij}$. 

\begin{prop}[See \cite{c1}, Proposition 1]\label{p:1}
Assume that $\mathcal{M}_\lambda$ in \eqref{e:3} satisfies the following condition
\begin{equation}\label{e:4}
\forall i,j\in \{1,2,\ldots,n\} \ \mathrm{such \ that} \ i\not= j \ \exists k\in \{1,2,\ldots, m\} \ \mathrm{such
\ that} \ d_i(k)\not= d_j(k).
\end{equation}
Then $\mathcal{Y}(U(I + \lambda G), \mathcal{M}_\lambda)=0$ for $\lambda \not=0$ small enough, where $G$ is an anti--hermitian matrix whose diagonal is null. Furthermore, its off--diagonal entries are \begin{equation} \label{e:5}g_{ij}=\frac{1}{2}\sum_{k=1}^mf^*_{ij}(k),
u^*_iR_ku_j+f_{ij}(k)u^*_iR^H_ku_j,
\end{equation}
 where
\begin{equation} \label{e:6}
f_{ij}(k)=\frac{d_j(k)-d_i(k)}{{\overset{m}{\underset{l=1} \sum}}|d_j(l)-d_i(l)|^2} \in \mathbb{C}.
\end{equation}
\end{prop}

\eqref{e:4} is essential in order to have a meaningful expression for \eqref{e:6}. 

\begin{rem} The original version of \cite[Proposition 1]{c1} is stated with $U(I + \lambda G+o(\lambda))J$, where $J$ is the product of a permutation matrix with a diagonal matrix having only unit modulus  entries on its diagonal and $o(\lambda)$ is a polynomial in $\lambda$ of degree $\ge 2$. In Proposition \ref{p:1} we have specialized $J=I$ and have considered only small enough $\lambda\not=0$ so that $o(\lambda)$ has been omitted. As noted in \cite[p.3 of \S 1.2]{c1}, if $U(I + \lambda G+o(\lambda))$ minimizes $\mathcal{Y}$, then also $U(I + \lambda G+o(\lambda))J$ minimizes $\mathcal{Y}$. In particular, this is true in the linear case, that is, when $\lambda$ is small enough and $o(\lambda)$ is omitted. All these observations justify our reformulation of Proposition \ref{p:1}. 
\end{rem}

The rest of the paper is devoted to modify Proposition \ref{p:1}, involving a wider family of solutions, minimizing $\mathcal{Y}$. Then we will be able to decompose $\mathrm{SL}(n,\mathbb{C})$, and consequently, $\mathrm{PSL}(n,\mathbb{C})$, $\mathrm{SU}(n)$ and $\mathrm{U}(n)$, in a canonical way (analogously, for the real case).


\section{Main results}
From \cite[Chapter 3]{r},   each invertible matrix with coefficients in $\mathbb{C}$ (resp. $\mathbb{R}$) and determinant equal to 1 can be generated by
matrices of the form $K(1,a)=I+ae_{ij}$, called \textit{transvection matrices}, where $a \in \mathbb{C}$ (resp. $a \in \mathbb{R}$). Note that $K(1,a)$ differs from  $I$ only in that
there is an $a$ in the $(i,j)$--th position, and, if $a=0$, then  $K(1,0)=I$. These matrices are important because  they generate $\mathrm{SL}(n,\mathbb{C})$ (resp. $\mathrm{SL}(n,\mathbb{R})$). More precisely, one can see that each $B \in  \mathrm{SL}(n,\mathbb{C})$ (resp. $\mathrm{SL}(n,\mathbb{R})$) can be decomposed uniquely as the product of finitely many $K(1,a)$, for suitable $a$. This is the well--known  \textit{rational canonical form of a special linear matrix} (see \cite[pp. 73--76]{r}).

\begin{lem}[3.2.10, see \cite{r}]  $\mathrm{SL}(n,\mathbb{C})=\langle K(1,a) \ | \ a \in \mathbb{C}\rangle$  and $\mathrm{SL}(n,\mathbb{R})=\langle K(1,a) \ | \ a \in \mathbb{R}\rangle$  for $n\ge2$.\end{lem}

Separately we do the following observation.

\begin{rem}\label{r:1}
$\mathrm{off}(K(1,a))=|a|^2$. If  $a=0$, then $\mathrm{off}(K(1,0))=\mathrm{off}(I)=0$.
The matrix $K(d_i(k),a)=d_i(k)+ae_{ij}$ is obtained by $I+ae_{ij}$, replacing the 1s on the principal diagonal with the
$d_1(k),\ldots,d_n(k)$, corresponding to $D_k$. If $d_i(k)=1$ for all $i=1,2,\ldots,n$, then
$K(d_i(k),a)=K(1,a)$. Then we can consider the family of matrices
\begin{equation}\label{e:7}
\mathcal{M}_{a,\lambda}= (U \cdot K(d_i(1),a) \cdot U^H + \lambda R_1, \ldots, U \cdot K(d_i(2),a) \cdot U^H + \lambda R_2,\ldots\end{equation}\[\ldots,
U \cdot K(d_i(m),a) \cdot U^H + \lambda R_m),\]  so that $\mathcal{M}_{0,\lambda}=\mathcal{M}_\lambda$ and $\mathcal{M}_{0,0}=\mathcal{M}_0$. Consequently, we can expect a significant generalization of the results in
\cite{a4,c2,c1}, once \eqref{e:3} is replaced by \eqref{e:7}.
\end{rem}

Now \cite[Lemma 1]{c1} may be reformulated by using algebraic methods and looking at $\mathcal{Y}$ from this new point of view. With the above notations, we  consider the linear maps
\begin{equation}\label{e:9}\gamma_{ijk}: V \in \mathrm{U}(n)\mapsto \gamma_{ijk}(V)=e^*_iV^HM^H_kVe_j\in
\mathbb{C}\end{equation}and\begin{equation}\label{e:9bis}   \ T_{ijk}: V \in \mathrm{U}(n) \mapsto
T_{ijk}(V)=e_ie^*_jV^HM^H_kV-V^HM_k^HVe_ie^*_j\in \mathrm{GL}(n,\mathbb{C}).
\end{equation}
Immediately, we deduce 
\begin{equation}\label{minimization}
\mathcal{Y}(V,\mathcal{M})= \sum^m_{k=1} \ \ {\underset{i\not=j}{\underset{i,j\in
\{1,\ldots,n\}}\sum}}|\gamma_{ijk}|^2.
\end{equation}
Now we  consider the linear map
\begin{equation}\label{e:8}
\mathcal{S}: (V,\mathcal{M})\in \mathrm{U}(n) \times \mathrm{GL}(n,\mathbb{C})^m \mapsto 
\mathcal{S}(V,\mathcal{M})=\sum^m_{k=1} \ \ {\underset{i\not=j}{\underset{i,j\in
\{1,\ldots,n\}}\sum}}\gamma^*_{ijk} \ T_{ijk}\in \mathrm{GL}(n,\mathbb{C}),\end{equation}
where the usual calculus rules $(A^*)^*=A$ and $(AB)^*=B^*A^*$ in $\mathrm{GL}(n,\mathbb{C})$ and  the equality $V^H=\overline{V}^*$, true by $V \in \mathrm{U}(n)$, imply
\begin{equation}\gamma^*_{ijk}(V)=(e^*_iV^HM^H_kVe_j)^*=e^*_jV^* (M^H_k)^* (V^H)^*e_i=
e^*_jV^* (M^H_k)^* \overline{V}e_i.\end{equation} 
In particular, if $(V,\mathcal{M})\in (\mathrm{U}(n) \cap \mathrm{GL}(n,\mathbb{R})) \times \mathrm{GL}(n,\mathbb{R})^m$, then $\overline{V}=V$, $V^H=V^*$, $(M^H_k)^*=\overline{M_k}=M_k$ and then
\begin{equation}\label{extra}\gamma^*_{ijk}(V)=e^*_jV^HM_k Ve_i.\end{equation} 
Once we replace  \eqref{extra} and \eqref{e:9bis} in \eqref{e:8}, we find exactly \cite[Equation 12]{c1}.
We may do similar considerations for the linear map
\begin{equation}\label{e:10}
\mathcal{S}^H: (V,\mathcal{M})\in \mathrm{U}(n) \times \mathrm{GL}(n,\mathbb{C})^m \mapsto 
\mathcal{S}^H(V,\mathcal{M})=\left( \sum^m_{k=1} \ \ {\underset{i\not=j}{\underset{i,j\in
\{1,\ldots,n\}}\sum}}\gamma^*_{ijk} \ T_{ijk}\right)^H\end{equation}
\[=\sum^m_{k=1} \ \ {\underset{i\not=j}{\underset{i,j\in \{1,\ldots,n\}}\sum}}(\gamma^*_{ijk})^H \ T^H_{ijk}\in \mathrm{GL}(n,\mathbb{C}).\]

A short proof of \cite[Lemma 1]{c1} is shown below.

\begin{lem}\label{l:1}  If $(V,\mathcal{M})\in (\mathrm{U}(n) \cap \mathrm{GL}(n,\mathbb{R})) \times \mathrm{GL}(n,\mathbb{R})^m$, then $\mathcal{S}$ is hermitian.
\end{lem}

\begin{proof} Since \eqref{extra} is satisfied, $(\gamma^*_{ijk})^H=\gamma^*_{ijk}$ and $T_{ijk}(V)^H=T_{ijk}(V)$ so we conclude
\begin{equation}\label{e:11}
\mathcal{S}(V,\mathcal{M})=\sum^m_{k=1} \ \ {\underset{i\not=j}{\underset{i,j\in
\{1,\ldots,n\}}\sum}}\gamma^*_{ijk} \ T_{ijk}=\sum^m_{k=1} \ \ {\underset{i\not=j}{\underset{i,j\in
\{1,\ldots,n\}}\sum}}\gamma_{ijk} \ T^H_{ijk}=\mathcal{S}^H(V,\mathcal{M}).\end{equation}
\end{proof}

Now we will generalize \cite[Lemmas 2, 3]{c1}, by looking at $\mathcal{M}_{a,\lambda}$ instead of $\mathcal{M}_\lambda$.
The condition $\mathcal{S}(V,\mathcal{M}_{a,\lambda})=\mathcal{S}^H(V,\mathcal{M}_{a,\lambda})$, described by Lemma \ref{l:1}, has a physical meaning. It is a stationary equation, which allow us to study the perturbations around the solution $(U, \mathcal{M}_{a,\lambda})\in (\mathrm{U}(n) \cap \mathrm{GL}(n,\mathbb{R})) \times \mathrm{GL}(n,\mathbb{R})^m$, minimizing $\mathcal{Y}$. A more detailed analysis can be done stopping at the terms of the first order in $\lambda\not=0$ and looking for an anti--hermitian $L\in \mathrm{GL}(n,\mathbb{C})$ such that, for small enough $\lambda\not=0$ and $a\not=0$, 
\begin{equation}\label{e:19}
\mathcal{S}(U(I+ \lambda L),\mathcal{M}_{a,\lambda})=\mathcal{S}^H(U(I+ \lambda L), \mathcal{M}_{a,\lambda}).\end{equation} 
The problem can be more conveniently centered at $U$, noting that
\begin{equation}\label{e:20}
\mathcal{S}(U (I+\lambda L),\mathcal{M}_{a,\lambda})=\mathcal{S}(I+\lambda L, \mathcal{N}_{a,\lambda}),
\end{equation}
where
\begin{equation}\label{n}
\mathcal{N}_{a,\lambda}=(K(d_i(1),a)+\lambda U^HR_1U, \ \  K(d_i(2),a)+\lambda U^HR_2U, \  \ldots \ldots \end{equation}
\[\ldots \ldots, K(d_i(k),a)+\lambda U^HR_kU, \  \ldots, \ K(d_i(m),a)+\lambda U^HR_mU).\]
\eqref{e:20} is exactly \cite[Equation 14]{c1}, once  $a=0$ and we stop at the terms of the first order in $\lambda$.

\begin{cor}\label{intermediate} From \eqref{e:20} and Lemma \ref{l:1}, for small enough $\lambda\not =0$ and $a \not=0$, \eqref{e:19} and $\mathcal{S}(I+ \lambda L,\mathcal{N}_{a,\lambda})=\mathcal{S}^H(I+ \lambda L, \mathcal{N}_{a,\lambda})$ have the same solutions.
\end{cor}

The following result provides more details.

\begin{lem}\label{l:3}Assume $(I+\lambda L,\mathcal{N}_{a,\lambda}) \in (\mathrm{U}(n) \cap \mathrm{GL}(n, \mathbb{R}) \times  \mathrm{GL}(n,\mathbb{R})^m$. Stopping at the terms of the first order in $\lambda\not=0$, the entries of $\mathcal{S}(I+\lambda L,\mathcal{N}_{a,\lambda})$, for small enough $\lambda\not=0$ and $a\not=0$, are
\begin{equation}\label{entries}s_{ij}(I+ \lambda L,\mathcal{N}_{a,\lambda})=\sum^m_{k=1} \delta_{ijk} \alpha_{ijk}+ (\delta_{ijk} \beta_{ijk} +\epsilon_{ijk} \alpha_{ijk}) \lambda,
\end{equation}
where $\alpha_{ijk}, \beta_{ijk}, \delta_{ijk}, \epsilon_{ijk}$ are linear maps depending only on $K(d_i(k),a)$ and $L$.
\end{lem}

\begin{proof}
From \eqref{e:9bis},
\begin{equation}\label{e:22} 
T_{ijk}(I+\lambda L)=e_ie^*_j(I+\lambda L)^H\cdot K(d_i(k),a)^H \cdot (I+\lambda L)-
\end{equation}
\[(I+\lambda L)^H \cdot K(d_i(k),a)^H \cdot (I+\lambda L) e_ie^*_j\]
\[=e_ie^*_j(I+\lambda L^H) \cdot K(d_i(k),a)^H \cdot (I+\lambda L)- (I+\lambda L^H) \cdot K(d_i(k),a)^H \cdot (I+\lambda L) e_ie^*_j\]
\[=(e_ie^*_jI+e_ie^*_j\lambda L^H) \cdot (K(d_i(k),a)^H \cdot I+K(d_i(k),a)^H \cdot \lambda L)\]\[-(I \cdot K(d_i(k),a)^H + \lambda L^H \cdot K(d_i(k),a)^H) \cdot (I\cdot e_ie^*_j+\lambda L \cdot e_ie^*_j) \]
\[=e_ie^*_jK(d_i(k),a)^H+\lambda e_ie^*_jK(d_i(k),a)^H L+\lambda e_ie^*_jL^HK(d_i(k),a)^H\]
\[+\lambda^2 e_ie^*_j L^H K(d_i(k),a)^H L - K(d_i(k),a)^H e_ie^*_j-\lambda K(d_i(k),a)^H L e_i e^*_j\]
\[ - \lambda L^H K(d_i(k),a)^H e_ie^*_j-\lambda^2 L^H K(d_i(k),a)^H L e_ie^*_j\]
Since we stop at the terms of the first order in $\lambda$, we  avoid $\lambda^2 e_ie^*_j L^H K(d_i(k),a)^H L$ and $\lambda^2 L^H K(d_i(k),a)^H L e_ie^*_j$, then
\[T_{ijk} (I+\lambda L)=e_ie^*_jK(d_i(k),a)^H+\lambda e_ie^*_jK(d_i(k),a)^H L+\lambda e_ie^*_jL^HK(d_i(k),a)^H\]
\[- K(d_i(k),a)^H e_ie^*_j-\lambda K(d_i(k),a)^H L e_i e^*_j- \lambda L^H K(d_i(k),a)^H e_ie^*_j.\]
Expanding with respect to $\lambda$, we put
\begin{equation}\label{e:23}\alpha_{ijk}=e_ie^*_jK(d_i(k),a)^H -K(d_i(k),a)^H e_i e^*_j\end{equation} 
and
\begin{equation}\label{e:24}\beta_{ijk}=(e_ie^*_jK(d_i(k),a)^H L - L^H K(d_i(k),a)^H e_i e^*_j)\end{equation}
\[+(e_ie^*_jL^HK(d_i(k),a)^H -K(d_i(k),a)^H Le_i e^*_j).\]
so that
\begin{equation}\label{e:25}
T_{ijk}(I+\lambda L)=\alpha_{ijk}+\beta_{ijk} \lambda. \end{equation}
On another hand, \eqref{extra} becomes
\begin{equation}\label{e:25}
\gamma^*_{ijk}(I+\lambda L)=e^*_i (I+\lambda L)^H  K(d_i(k),a) (I+\lambda L) e_j\end{equation}
\[=e^*_i (I+\lambda L^H)  K(d_i(k),a) (I+\lambda L) e_j\]
\[=(e^*_i K(d_i(k),a) +\lambda e^*_i L^HK(d_i(k),a)) \cdot (Ie_j+\lambda Le_j)\]
\[=e^*_i K(d_i(k),a) e_j + \lambda e^*_i K(d_i(k),a)L e_j+\lambda e^*_i L^H K(d_i(k),a) e_j+ \lambda^2 e^*_i L^H K(d_i(k),a) L e_j.\]
Since we stop at the terms of the first order in $\lambda$, we avoid $\lambda^2 e^*_i L^H K(d_i(k),a) L e_j$, then
\[=e^*_i K(d_i(k),a) e_j + \lambda e^*_i K(d_i(k),a)L e_j+\lambda e^*_i L^H K(d_i(k),a) e_j.\]
Expanding with respect to $\lambda$, we put
\begin{equation}\delta_{ijk} = e^*_i K(d_i(k),a) e_j
\end{equation}
and
\begin{equation}\epsilon_{ijk} = e^*_i K(d_i(k),a) L e_j + e^*_i L^H K(d_i(k),a) e_j
\end{equation}
so  that
\begin{equation}\label{e:30}
\gamma^*_{ijk}=\delta_{ijk}+\lambda \epsilon_{ijk}.
\end{equation}
Collecting all we have  found, 
\begin{equation}\label{e:31}
s_{ij}(I+ \lambda L,\mathcal{N}_{a,\lambda})=\sum^m_{k=1} \ \ {\underset{i\not= j}{\underset{i,j\{1,\ldots, n\}} \sum}}\gamma^*_{ijk} \ T_{ijk} \end{equation} 
\[=\sum^m_{k=1} \ \ {\underset{i\not= j}{\underset{i,j\{1,\ldots, n\}} \sum}} (\delta_{ijk} + \epsilon_{ijk} \lambda) (\alpha_{ijk}+\beta_{ijk} \lambda).\]
Now $(\delta_{ijk} + \epsilon_{ijk} \lambda) (\alpha_{ijk}+\beta_{ijk} \lambda)=\delta_{ijk} \alpha_{ijk} + (\delta_{ijk}\beta_{ijk}+\epsilon_{ijk}\alpha_{ijk}) \lambda + (\epsilon_{ijk}\beta_{ijk}) \lambda^2=\delta_{ijk} \alpha_{ijk} + (\delta_{ijk}\beta_{ijk}+\epsilon_{ijk}\alpha_{ijk}) \lambda$, since we avoided the term $(\epsilon_{ijk}\beta_{ijk}) \lambda^2$, as done until now. The result follows.
\end{proof}

Now we have all that is necessary to state our main result.

\begin{thm}\label{t:1}
Assume $(I+\lambda L,\mathcal{N}_{a,\lambda}) \in (\mathrm{U}(n) \cap \mathrm{GL}(n, \mathbb{R}) \times  \mathrm{GL}(n,\mathbb{R})^m$.
Then, for small enough $\lambda \not=0$ and $a\not=0$, we get around $U \in \mathrm{U}(n) \cap \mathrm{GL}(n, \mathbb{R})$  
\begin{equation}\mathcal{Y}(I + \lambda L, \mathcal{N}_{a,\lambda})=\sum^m_{k=1} \ \ {\underset{i\not=j}{\underset{i,j\in
\{1,\ldots,n\}}\sum}}|\gamma_{ijk}(I + \lambda L)|^2,\end{equation} where $L$ is an anti--hermitian matrix such that the stationary condition $\mathcal{S}(I+ \lambda L,\mathcal{N}_{a,\lambda})=\mathcal{S}^H(I+ \lambda L, \mathcal{N}_{a,\lambda})$ is satisfied with entries 
$s_{ij}(I+ \lambda L,\mathcal{N}_{a,\lambda})$ as in \eqref{entries}. 
\end{thm} 

\begin{proof}  The result follows from Corollary \ref{intermediate} and Lemma \ref{l:3}. \end{proof}

The reader may note that the cases $\lambda=a=0$ and $\lambda\not=0$ with $a=0$ have been already described in \cite{a4,c2,c1} and recalled above. The  choice of $\mathcal{M}_{a,\lambda}$, and consequently of $\mathcal{N}_{a,\lambda}$, is justified by the next result. Preliminarly, we note that
\begin{equation}
\mathcal{N}_{a,\lambda}=(K(d_i(1),a), \ldots , K(d_i(m),a))+\lambda(U^HR_1U, \ldots , U^HR_mU)=\mathcal{C}_a+\lambda \mathcal{D},
\end{equation}
where $\mathcal{C}_a \in \mathrm{SL}(n,\mathbb{C})^m$ and $\mathcal{D}\in \mathrm{GL}(n,\mathbb{C})^m$ in general.

\begin{cor}\label{conclusion} In the hypotheses of Theorem \ref{t:1}, we may draw the same conclusions, replacing $\mathcal{N}_{a,\lambda}$ with an arbitrary $\mathcal{N} =(\mathcal{C}_{a_1}+\lambda_1\mathcal{D}_1) \cdot \ldots  \cdot (\mathcal{C}_{a_r}+\lambda_r\mathcal{D}_r)$, for small enough $a_1,\ldots, a_r$ and $\lambda_1, \ldots, \lambda_r$, where $r\geq1$.
\end{cor}

\begin{proof}
Lemma \ref{l:1} mention the rational canonic form of a special linear matrix $\mathcal{C} \in \mathrm{SL}(n,\mathbb{C})$ so we have  a  unique decomposition of $\mathcal{C}$ as product of $\mathcal{C}_{a_1}, \ldots, \mathcal{C}_{a_r}$. Consequently, we have a unique decomposition of $\mathcal{N}$ as claimed.  Theorem \ref{t:1} is true for each factor and  the linearity of $\mathcal{Y}$ and $\mathcal{S}$ allows us to end.
\end{proof}

An intriguing open problem is stated at this point. Some times we can consider families of matrices of $\mathrm{SL}(n,\mathbb{C})$  which are dense in $\mathrm{GL}(n,\mathbb{C})$. When this happens, in principle, we may approximate any matrix of $\mathrm{GL}(n,\mathbb{C})$ with a suitable sequence of matrices of $\mathrm{SL}(n,\mathbb{C})$, for which we have a description by Corollary \ref{conclusion}. Unfortunately, we were not able to provide an efficient algorithm of  approximation of a suitable family of matrices of $\mathrm{SL}(n,\mathbb{C})$, dense in $\mathrm{GL}(n,\mathbb{C})$.


\end{document}